\documentclass[10pt]{article}
\usepackage{amsfonts, amssymb, amsmath, amscd}
\usepackage{graphicx}
\usepackage{latexsym, longtable, tabularx, array}
\usepackage{epsfig, psfrag}
\usepackage[usenames]{color}

\newtheorem{thm}{Theorem}

\newtheorem{lemma}{Lemma}
\newtheorem{cor}[thm]{Corollary}

\newtheorem{rmk}{Remark}

\newcommand{\C}{\mathbb{C}}
\newcommand{\cl}[1]{\ensuremath{\overline{ #1}}}
\newcommand{\kb}[1]{\ensuremath{\langle #1 \rangle}}
\newcommand{\D}{\Delta^2}
\newcommand{\pf}{{\em Proof: \quad }}
\newcommand{\done}{\hfill $\blacksquare$}
\newcommand{\myone}{\ensuremath{{\bf 1}}}
\newcommand{\ZA}{{\mathbb{Z}}[A^{\pm1}] }
\newcommand{\SL}{\mathfrak s\mathfrak l}

\hoffset=- 1 cm
\begin{document}
\title{On the tail of Jones polynomials\\ of closed braids with a full twist}
\author{
Abhijit  Champanerkar and  Ilya Kofman 
\footnote{The authors gratefully acknowledge support by the NSF and PSC-CUNY.} \\
{\em {\small Department of Mathematics, College of Staten Island,
 City University of New York}}} 
\maketitle
\begin{abstract}\noindent
  For a closed $n$--braid with a full positive twist and with $\ell$
  negative crossings, $0\leq \ell \leq n$, we determine the first
  $n-\ell+1$ terms of the Jones polynomial $V_L(t)$.  We show that
  $V_L(t)$ satisfies a braid index constraint, which is a gap of
  length at least $n-\ell$ between the first two nonzero coefficients
  of $(1-t^2) V_L(t)$. For a closed positive $n$--braid with a full
  positive twist, we extend our results to the colored Jones
  polynomials.  For $N>n-1$, we determine the first $n(N-1)+1$ terms of
  the normalized $N$--th colored Jones polynomial.
\end{abstract}

\section{Introduction}

The {\em tail} (resp. {\em head}) of a polynomial will denote the
sequence of its lowest (resp. highest) degree terms, up to some
specified length.  In this note, we precisely determine the tail of
the Jones polynomial for a closed $n$--braid with a full positive
twist, and with up to $n$ negative crossings.  We also precisely
determine the tails of the colored Jones polynomials for a closed
positive $n$--braid with a full positive twist.

It is natural to consider quantum and geometric invariants of links
that are closed braids with a full twist.  For example, Lorenz links,
all of which are closed positive braids with a full twist, dominate
the census of the simplest hyperbolic knots, and their Jones
polynomials are relatively simple \cite{BK, CKP}.  The full twist
arises as $\pm 1$ Dehn surgery on the braid axis, considered as an
augmented unknot in $S^3$.  Hence, adding full twists is a natural
geometric operation on links.  On the other hand, the full twist is in
the center of the braid group, so its image in any irreducible
representation is a scalar.  All known closed formulas for Jones
polynomials of infinite link families essentially rely on this fact.

For any closed braid, we showed in \cite{mjp} that after sufficiently
many full twists on a subset of strands, the coefficient vector for
any colored Jones polynomial decomposes into fixed blocks, separated
by blocks of zeros that increase by a constant length for every twist.
So once the non-zero blocks separate, they simply move apart unchanged
with every additional full twist.  In Theorem \ref{mainthm1} below, we
completely determine the first block for full twists on all strands.
In this case, the first block separates after only one full twist.

Dasbach and Lin \cite{DL2} showed that for alternating knots, and more
generally $A$--adequate knots, the first three coefficients in the
tail of the normalized $N$-th colored Jones polynomials are
independent of the color for $N\geq 3$.  In Corollary \ref{tailcor}
below, for $N>n-1$, we determine the tail of length $n(N-1)+1$ for the
normalized $N$--th colored Jones polynomial of closed positive braids
with a full twist. In fact, Theorem \ref{mainthm2} implies that given
$M \geq 2$, for all colors $N\geq M$, the coefficients of the tail of
length $M$ stabilize up to sign.

Dasbach and Lin also showed that the second coefficients of the head
and tail together provide a linear bound for the hyperbolic volume of
alternating knots.  In contrast, the coefficients of the tail of
length $N$ as in Theorem \ref{mainthm2} for closed positive braids
with a full twist are all $\{-1,0,1\}$.  Moreover, the tail of length
$n(N-1)$ as in Corollary \ref{tailcor}, also has coefficients only
$\{-1,0,1\}$.  These coefficients and the dependence on the braid
index indicates that, for this class of knots, these tails by
themselves are unrelated to the hyperbolic volume of the closed braid.
For example, $3$--braids can have unbounded hyperbolic volume, and the
positive twisted torus knots $T(p,q,2,2k)$ have bounded hyperbolic
volume but unbounded braid index \cite{ttk}.

To state our main results, we adopt the following standard
convention.  
Let $V_L(t)$ denote the Jones polynomial, such that
$$ t^{-1}V_{L_+} - tV_{L_-} = \left(t^{1/2}- t^{-1/2}\right)V_{L_0}
\quad {\rm and} \quad V_{\bigcirc}(t)=1. $$

\begin{thm}\label{mainthm1}
  Let $\beta'$ be a $n$--braid of length $c$ with $\ell$ negative crossings  with 
$0\leq \ell \leq n$ and $\beta=\D_n
  \,\beta'$, where $\D_n$ is the positive full twist in the braid
  group $B_n$. Then
\begin{align*}
 V_{\cl{\beta}}(t) & = (-1)^{n+c+1}\;\,t^{\frac{(n-1)^2+c-2\ell}{2}}
\left(\frac{1 + t^{n-\ell+1}\; p(\cl{\beta};\, t)}{1-t^2} \right)\\
& =  (-1)^{n+c+1}\;\,t^{\frac{(n-1)^2+c-2\ell}{2}}
\left(\sum_{i=0}^{[(n-\ell)/2]} t^{2i} + t^{n-\ell +1}q(\cl{\beta};\, t)\right) 
\end{align*}
where $p(\cl{\beta};\, t)$ and $q(\cl{\beta};\, t)$ are polynomials in $t$. 
\end{thm}

The latter expression gives the tail of $V_{\cl{\beta}}(t)$ of length
$n-\ell+1$.  An interesting consequence is that the Jones polynomial
satisfies a braid index constraint, which is a gap of length at least
$n-\ell$ between the first two nonzero coefficients of $(1-t^2)\cdot
V_{\cl{\beta}}(t)$.  If $\ell=0$ in Theorem \ref{mainthm1}, then
$\beta$ is a positive $n$--braid with a full twist. In this case, the
MFW inequality \cite{Morton:86,FW:87} is sharp so the braid index of
$\cl{\beta}$ is $n$. However, the gap between the first two nonzero
coefficients of $(1-t^2)V_{\cl{\beta}}(t)$ can be more than $n$.  For
example, if $\beta'=\sigma_2^2 \sigma_1 \in B_4$ then $\beta=\D_4
\beta'$ and
$$V_{\cl{\beta}}(t)= t^6 +t^8 +t^{10} +t^{12} \ \ \implies\ \ (1-t^2) V_{\cl{\beta}}(t) = t^6 -t^{14}.$$

Another consequence of Theorem \ref{mainthm1} is related to a
conjecture of V. Jones \cite{Jones}, which remains open in general:
The writhe $w(\beta)$, which is the algebraic crossing number of
$\beta$, is a topological invariant of $\cl{\beta}$ whenever $n$ is
the minimal braid index of $\cl{\beta}$ (see \cite{Stoimenow2002}).
When the MFW inequality is sharp, the Jones conjecture is known to
hold, which is the case for positive braids, $\ell=0$.  If $\ell=1$,
$\beta$ is conjugate to a positive braid.  For $\ell>1$, although we
do not know when the MFW inequality is sharp, we can prove the Jones
conjecture:
\begin{cor}
For $\beta$ as in Theorem \ref{mainthm1}, 
$$ 2\min\deg(V_{\cl{\beta}}(t))= w(\beta)-n+1. $$ 
Thus, if $n$ is the minimal braid index of $\cl{\beta}$, then
$w(\beta)$ is a topological invariant of $\cl{\beta}$.
\end{cor}
%\pf
{\em Proof: } 
$2\min\deg(V_{\cl{\beta}}(t))=(n-1)^2+c-2\ell = n(n-1)+c-2\ell-n+1 = w(\beta)-n+1$.
\done

Let
$J_N(L; t)$ be the colored Jones polynomial of $L$, colored by the
$N$-dimensional irreducible representation of $\SL_2(\C)$, with the
normalization 
$$J_2(L;\,t)=(t^{1/2}+t^{-1/2})V_L(t) \ \mathrm{ and} \ \ 
J_N(\bigcirc;\,t) = \frac{t^{N/2}-t^{-N/2}}{t^{1/2}-t^{-1/2}} =[N].
$$ 
The colored Jones polynomials are weighted sums of Jones polynomials
of cablings, and the following formula is given in \cite{KM}.  Let
$L^{(r)}$ be the $0$-framed $r$-cable of $L$; i.e., if $L$ is
$0$-framed, then $L^{(r)}$ is the link obtained by replacing $L$ with
$r$ parallel copies.  (See below for a formula modified for other
framings.)
\begin{equation}\label{colJP}
J_{N+1}(L;\,t) = \sum^{[N/2]}_{j=0}(-1)^j\binom{N-j}{j}J_2({L^{(N-2j)}};\,t)
\end{equation}
The normalized colored Jones polynomial $J'_N(L;\,t)$ is defined by 
$\displaystyle J'_N(L;\,t)=\frac{J_N(L;\,t)}{[N]}$ and $\displaystyle J'_N(\bigcirc;\, t)=1.$

\begin{thm}
\label{mainthm2}
  Let $\beta'$ be a positive $n$--braid of length $c$ and $\beta=\D_n
  \,\beta'$, where $\D_n$ is the positive full twist in the braid
  group $B_n$. Then
\begin{align*}
 J_{N+1}(\cl{\beta};\,t) &= (-1)^{N(n+c+1)}t^{\frac{N((n-1)^2+c-1)}{2}} 
\left(\frac{1+t^{nN+1}p_N(\cl{\beta};\, t)}{1-t}\right) \\
 &= (-1)^{N(n+c+1)}t^{\frac{N((n-1)^2+c-1)}{2}} 
\left(\sum_{i=0}^{nN} t^i + t^{nN+1} q_N(\cl{\beta};\, t)  \right)
\end{align*}

\begin{align*}
J'_{N+1}(\cl{\beta};\,t) &= (-1)^{N(n+c+1)}t^{\frac{N((n-1)^2+c)}{2}} 
\left(\frac{1+t^{nN+1}p_N(\cl{\beta};\, t) }{1-t^{N+1}}\right)\\
&= (-1)^{N(n+c+1)}t^{\frac{N((n-1)^2+c)}{2}} 
\left(\sum_{i=0}^{n-1}t^{i(N+1)} + t^{nN+1}q_N(\cl{\beta};\, t)\right)
\end{align*}
where $p_N(\cl{\beta};\, t)$ and $q_N(\cl{\beta};\, t)$ are polynomials in $t$. 
\end{thm}

\begin{cor}\label{tailcor}
If $\beta'$ is a positive $n$--braid of length $c$ and $\beta=\D_n
\,\beta'$, then for $N>n-2$,
$$ J'_{N+1}(\cl{\beta};\,t) = (-1)^{N(n+c+1)}t^{\frac{N((n-1)^2+c)}{2}} \left( 
\sum_{i=0}^{n-1}t^{i(N+1)} \ + \ \{{\rm terms \ of \ degree } \geq nN+1 \}\right). $$
\end{cor}
\pf
In Theorem \ref{mainthm2}, $\sum_{i=0}^{n-1}t^{i(N+1)}$ and $t^{nN+1} q_N(\cl{\beta};\, t)$ can overlap  
only when\\ $(n-1)(N+1) \geq nN+1$, i.e. when $N \leq n-2$. \done

For $N>n-1$, Corollary \ref{tailcor} determines the tail of length
$n(N-1)+1$ for the normalized $N$--th colored Jones polynomial.

The Jones Slope Conjecture \cite{stavros2010} relates boundary slopes of a knot $K$ to 
the growth rate of the lowest and highest degrees of $J'_N(K;\,t)$. 

\begin{cor}
\label{slopecor}
Let $j^*(N)$ denote the lowest degree of $J'_N(\cl{\beta};\,t)$ for
$\beta$  as in Theorem \ref{mainthm2}. Then 
$$\lim_{N\to \infty} \frac{4j^*(N)}{N^2}=0$$ 
This proves part of the Jones Slope Conjecture for knots $\cl{\beta}$.
\end{cor}
\pf 
By Corollary \ref{tailcor}, $j^*(N)=(N-1)((n-1)^2+c)/2$ for $N>n-2$ which 
proves the limit. 
Since $\beta$ is a positive braid, $\cl{\beta}$ is fibered and the
fiber has slope $0$, thus proving part of the Jones Slope Conjecture. 
\done

 Since
$\cl{\beta}$ is $A$-adequate, Corollary \ref{slopecor} also follows from Example 9 in
\cite{KFP2010}.

 \begin{rmk}\rm The AJ Conjecture \cite{stavrosAJC} relates the recurrence 
   properties of $J'_N(K;t)$ to the A-polynomial of a knot
   $K$.  The terms of $J'_N(K;t)$ in Corollary \ref{tailcor} have
   exponents which are linear in the color $N$, so these terms give no
   information about the A-polynomial via the AJ-Conjecture.
 \end{rmk}

\section{Proof of Theorem \ref{mainthm1}}

Generalizing the well-known formula for torus knots, the Jones
polynomial of any torus link $T(p,q)$ is given by the following sum,
with $d=\gcd(p,q)$ \cite{ILR}.

\begin{equation}\label{JPtorus}
V_{T(p,q)}(t) = (-1)^{d+1}
\frac{t^{(p-1)(q-1)/2}}{1-t^2}\sum_{i=0}^d{d\choose i}t^{\frac{p}{d}(1+\frac{q}{d}i)(d-i)}\left(t^{\frac{q}{d}(d-i)}-t^{1+\frac{q}{d}i}\right)
\end{equation}

For $n=2$, the claim follows from (\ref{JPtorus}), which in this case simplifies to 
$$V_{T(2,q)}(t)=(-1)^{q+1}\frac{ t^{(q-1)/2}}{(1-t^2)}(1-t^3+(-1)^q(t^{1+q}-t^{2+q})).$$

Henceforth, let $n>2$.
The Jones polynomial $V_L(t)$ is obtained from the Kauffman bracket $\kb{L}$ by
substituting $t = A^{-4}$ and multiplying by $(-A^{3})^{-w}$ to adjust for the writhe 
$w$ of $L$. We will show that the right-most part of
the coefficient vector of the Kauffman bracket 
$\vec{c}=\{c_i\, |\, \kb{L}=A^{*}\sum c_i\, A^{4i}\, \}$
has the following form:
\begin{align*}
* \;\underbrace{ 0\; 1\; \ldots 0\; 1\;}_{n+1} & \ {\rm if} \ n \ {\rm is\ odd}\  & \qquad &
* \;\underbrace{-1\; 0\; -1\; \ldots 0\; -1\;}_{n+1} \ {\rm if} \ n \ {\rm is\ even}\ 
 \end{align*}
We will call this part of the Kauffman bracket the \emph{gap block}. 
Multiplying by $1-t^2 =1-A^{-8}$, the coefficient vector changes as follows: 

\begin{align*}
\ &\qquad \fbox{Case: $n$ odd} & \fbox{Case: $n$ even}\qquad \\
{\rm coefficients\ of}\left((-1)^{n+1}\kb{L}\right) \ &\quad \ *\ *\ *\ 0\ 1 \ldots 0\ 1\ 0\ 1\  & \  *\ *\ *\ 1\ 0\ 1 \ldots 0\ 1\ 0\ 1\  \\
{\rm coefficients\ of}\left((-1)^{n+1}A^{-8}\kb{L}\right)\ &\quad \ *\ 0\ \ 1\ 0\ 1 \ldots 0\ 1\ \ \,  & 1\ \; 0\  1\ 0\ 1 \ldots 0\ 1\qquad  \\
 \hline
{\rm coefficients\ of}\left((-1)^{n+1}(1-A^{-8})\kb{L}\right)\ &\quad \ *\ *\ *\ \underbrace{0\ 0 \ldots 0\  0\  0\ }_{n} 1\  & \ *\ *\ *\ \ \underbrace{0\ 0\, \ldots\, \ 0\  0\  0\ }_{n} 1\  
\end{align*}

When $p=q=n$ in (\ref{JPtorus}), we obtain the Jones polynomial of
$(n,n)$ torus link $=\cl{\D_n}$, which is the closure of the full twist in the braid group $B_n$. 

\begin{align}\label{JPtorus2}
 V_{\cl{\D_n}}(t) \ =\ & 
(-1)^{n+1} \frac{t^{\frac{1}{2}(n-1)^2}}{1-t^2}\sum_{i=0}^n{n\choose i}t^{(1+i)(n-i)}\left(t^{n-i}-t^{1+i}\right) \notag\\
=\ & (-1)^{n+1} \frac{t^{\frac{1}{2}(n-1)^2}}{1-t^2}\left(1+(n-2)t^{n+1}\ + \ \{{\rm higher\ order\ terms}\} \right) 
\end{align}

The Kauffman bracket of $\cl{\D_n}$ is obtained from (\ref{JPtorus2})
by substituing $t=A^{-4}$ and to adjust for the writhe, multiplying by $(-A)^{3n(n-1)}=A^{3n^2-3n}$.  The
first two terms of $(1-t^2)\cdot V(t)$ change as follows:
\begin{align*}
A^{3n^2-3n}\cdot A^{-2(n-1)^2}\left(1+(n-2)A^{-4(n+1)}\right) % =A^{3n^2-3n}\cdot(A^{-2n^2+4n-2}+(n-2)A^{-2n^2-6})
= A^{n^2+n-2}+(n-2) A^{n^2-3n-6}
\end{align*}
% become the following terms: $(n-2) A^{n^2-3n-6}$ and $A^{n^2+n-2}$.  
After dividing by $1-A^{-8}$ we obtain a sum that depends on the parity of $n$. We
see that the Kauffman bracket of the positive full twist on
$n$-strands has a gap block with top degree $n^2+n-2$, and its smallest non-zero coefficient has degree $n^2-3n\pm 2$, according to whether $n$ is even or odd.
However, when $n$ is odd, the zero ``term'' is the first one in the gap block, so
we say that the gap block has bottom degree $n^2-3n-2$; i.e.
\begin{equation}\label{KBtorus2}
(-1)^{n+1}\kb{\cl{\D_n}}=  q_2(A) + a_n\;A^{n^2-3n-2} + A^{n^2-3n+2} + \ldots + A^{n^2+n-2}
\end{equation}
where $q_2(A)$ is a Laurent polynomial with degree strictly less than $n^2-3n-2$ and 
$a_n =(1+(-1)^n)/2$.  
We will show that adding the positive braid $\beta'$ to $\D_n$ does
not affect the gap block.

The Temperley-Lieb algebra $TL_n$ is closely related to the Jones polynomial.
In the usual notation, $TL_n$ is the algebra over $\ZA$ with generators $\{\myone, e_1, e_2,\ldots,e_{n-1}\}$ and relations, with $\delta=-A^2-A^{-2}$,
$$ e_i^2=\delta e_i,\quad e_i e_{i\pm 1} e_i = e_i,\quad e_i e_j = e_j e_i \text{ if } |i-j|\geq 2 $$

As a free $\ZA$-module, $TL_n$ has a basis that consists of all
diagrams with no crossings and no closed curves, with dimension equal
to the Catalan number $C_n = \frac{1}{n+1} \binom{2n}{n}$. We will
refer to this particular basis as the {\em standard basis} $\{h_i\,
|\, i=0, \ldots, C_n-1 \, \}$ with $h_0 = \myone$.  Each $h_i$ can be
expressed as a product of distinct generators: $h_i=e_{j_1} \ldots
e_{j_{r}}$.

Let $A^{a-b}$ be the contribution from any
smoothing $s$ of $\beta'$ where $a$ and $b$ are the number of $A$ and $B$
smoothings of $\beta'$, respectively.  Passing to the representation of $\beta'$ in $TL_n$, 
$$\beta'=\sum\nolimits_s A^{a(s)-b(s)}\; \bigcirc^{|s|} \; h_{i(s)}$$ 
where $h_{i(s)}$ is the basis element obtained 
from a smoothing $s$, and $|s|$ is the number of loops in the smoothing 
of $\beta'$ (not the closure of $\beta'$).   

Let $c={\rm length}(\beta')$. 
For $ 0\leq \ell \leq n$, the (unique) state which gives $h_0$ has $\ell$ $B$--smoothings, for
which $a-b=(c-\ell)-\ell=c-2\ell$.  
We define $q_1(A)$ as follows:
\begin{align*}
(-1)^{n+1}q_1(A)&=\kb{\cl{\D_n \beta'}}-A^{c-2\ell}\kb{\cl{\D_n h_0}} \\
&= \sum_{s\ \text{with}\ i(s)\neq 0} A^{a(s)-b(s)} \kb{\bigcirc^{|s|}\;\cl{\D_n h_{i(s)}}}\\
&= \sum_{s\ \text{with}\ i(s)\neq 0} A^{a(s)-b(s)} \delta^{|s|} \kb{\;\cl{\D_n h_{i(s)}}}
\end{align*}

By Lemma \ref{maxpower}, which is proved below, for $i>0$ the highest power in
$\kb{\;\cl{\D_n h_{i}}}$ is $n^2-3n-4$. This implies the following:

\begin{lemma}
\label{maxpowerq}
If $0\leq\ell\leq n$, the degree of $q_1(A)$ is at most $c+n^2-3n-6+2\ell$.
\end{lemma}
\pf
First, suppose $\ell=0$.
We claim that to get $k$ loops in any smoothing of $\beta'$ we need at
least $(k + 1)$ $B$--smoothings.  Since $\beta'$ is a positive braid,
every $B$--smoothing adds at most one loop, but the first
$B$--smoothing does not result in any loops.  Hence, $(k+1)$ $B$--smoothings
(and the remaining $A$--smoothings) result in at most $k$ loops.  

It follows that $1 \leq k+1 \leq b$ and $ a \leq c - k -1$.
Hence, $a-b \leq c-2k-2$ and the highest power in $\delta^k$ is
$2k$. Thus, for $\ell=0$, the highest power in $A^{a(s)-b(s)} \delta^{|s|}$ is
$c-2k-2+2k = c-2$.  

If $\ell>0$, we claim that to get $k$ loops in any smoothing of
$\beta'$ we need at least $(k + 1 - \ell)$ $B$--smoothings.  As for a
positive braid, $(k+1)$ smoothings that produce a cup-cap give at most
$k$ loops. But now, some of these smoothings could be $A$--smoothings
at a negative crossing, so $(k+1-\ell)$ $B$--smoothings (and the
remaining $A$--smoothings) result in at most $k$ loops.

It follows that $b \geq k+1-\ell $ and $ a \leq c - k -1+\ell$.  Hence, $a-b \leq
c-2k-2+2\ell$.  Since the highest power in $\delta^k$ is $2k$, the
highest power in $A^{a(s)-b(s)} \delta^{|s|}$ is $c-2k-2+2\ell+2k =
c+2\ell-2$. 

By Lemma \ref{maxpower}, the highest power in
$\kb{\cl{\D_n h_{i(s)}}}$ is $n^2-3n-4$, so the degree of $q_1(A)$ is at most $c+n^2-3n-6+2\ell$. 
\done

We now return to the equation, $(-1)^{n+1}q_1(A)=\kb{\cl{\D_n \beta'}}-A^{c-2\ell}\kb{\cl{\D_n}}$.
By Lemma \ref{maxpowerq}, the degree of $q_1(A)$ is at most $\alpha_1=c+n^2-3n-6+2\ell$.
Hence by (\ref{KBtorus2}),
\begin{equation}\label{bracketofbeta}
(-1)^{n+1}\kb{\cl{\D_n \beta'}}=q_1(A) + A^{c-2\ell} \left( q_2(A) + a_n\;A^{n^2-3n-2}
+ A^{n^2-3n+2} + \ldots + A^{n^2+n-2} \right)
\end{equation}
The highest power of $A$ above is $n^2+n-2+c-2\ell$.
Thus, the gap block starts at the power ${\alpha_2=c+n^2-3n-2-2\ell}$. 
Comparing with the maximum power of $q_1(A)$, we see that $q_1(A)$ and the gap block can
overlap in at most $\ell$ coefficients:
$\alpha_1-\alpha_2=4(\ell-1)$.
% =c+n^2-3n-4+2\ell-2-(c+n^2-3n-2-2\ell)=4\ell-4=4(\ell-1)$.
After multiplying by $(1-t^2)$, we get a gap of $(n-\ell)$ zeros.

Let us compute the highest power of $A$ after 
adjusting for the writhe. The writhe $w=n(n-1)+(c-\ell)
-\ell=n^2-n+c-2\ell$. So after multiplying by $(-A^3)^{-w}$, the highest
power of $A$ is
\begin{align*}
(-1)^{n+1} A^{n^2+n-2+c-2\ell}(-1)^{n^2-n+c-2\ell}A^{-3(n^2-n+c-2\ell)}&=
(-1)^{n^2+c+1}A^{-2n^2+4n-2c+4\ell-2} \\
&= (-1)^{n+c+1}(A^{-4})^{\frac{(n-1)^2+c-2\ell}{2}}
\end{align*}
To get the Jones polynomial, we substitute $t=A^{-4}$ and multiply by
$(-A^3)^{-w}$. Hence we obtain the lowest power of $t$ to be 
$(-1)^{n+c+1}t^{\frac{(n-1)^2+c-2\ell}{2}}$.
This completes the proof of the first statement of Theorem \ref{mainthm1}.  

To obtain the tail (without denominators), we note that
$\displaystyle{\frac{1+tp(t)}{1-t^2} = 1+tq(t)}$, where $q(t)$ is a
polynomial.  For simplicity, let $\ell=0$, but the proof is same in
the other case.  Given the polynomial $p(t)$ from the proof above, we
obtain the polynomial $q(t)$ as follows:

\begin{align*}
n \ {\rm odd} \hspace*{1in} \frac{1 + t^{n+1}\; p(t)}{1-t^2} %& = \frac{1 -t^{n+1} + t^{n+1}(1+p(t))}{1-t^2}\\
&= \frac{1-t^{n+1}}{1-t^2} + \frac{t^{n+1}(1+p(t))}{1-t^2} \\
&= \sum_{i=0}^{[n/2]} t^{2i} + t^{n+1}q(t)
\end{align*} \\
\begin{align*}
n \ {\rm even} \hspace*{1in} \frac{1 + t^{n+1}\; p(t)}{1-t^2} %& = \frac{1 -t^{n} + t^{n}(1+t p(t))}{1-t^2}\\
&= \frac{1-t^{n}}{1-t^2} + \frac{t^{n}(1+t p(t))}{1-t^2} \\
&= \sum_{i=0}^{(n-2)/2} t^{2i} + t^{n}(1+t q(t)) \\
&= \sum_{i=0}^{n/2} t^{2i} + t^{n+1} q(t)
\end{align*}

This completes the proof of Theorem \ref{mainthm1}.  
\done

\begin{lemma}
\label{fulltwistevaluation}
Let $h_i$ be a standard basis element of $TL_n$. 
Let $k = \#$cups in $h_i$, and $m = \#$ through strands in $h_i$, so that $2k+m=n$.
If $H_i = h_i$ with its $m$ through strands given a full right twist, then
$\D_n h_i$ = $\D_n$ if $i=0$, and $\D_n h_i = A^{-6k} H_i$ if $i>0$.
\end{lemma}
\pf 
Let $B=\{ 1,2, \ldots,2n\}$ denote the positions of the strands.
Let ${\rm cap}(h)=\{ u_1,\ldots,u_{2k} \}$ denote indices of strands of $h$ which are caps; e.g., ${\rm cap}(e_i)=\{i,i+1\}$. 
Let ${\rm feet}(h)= \{v_1,\ldots,v_m \}$ denote the bottom indices of strands of $h$  which pass through. 
Note that ${\rm cap}(h) \sqcup {\rm feet}(h) =B$.

\begin{figure}
\begin{center}
\includegraphics[height=2in]{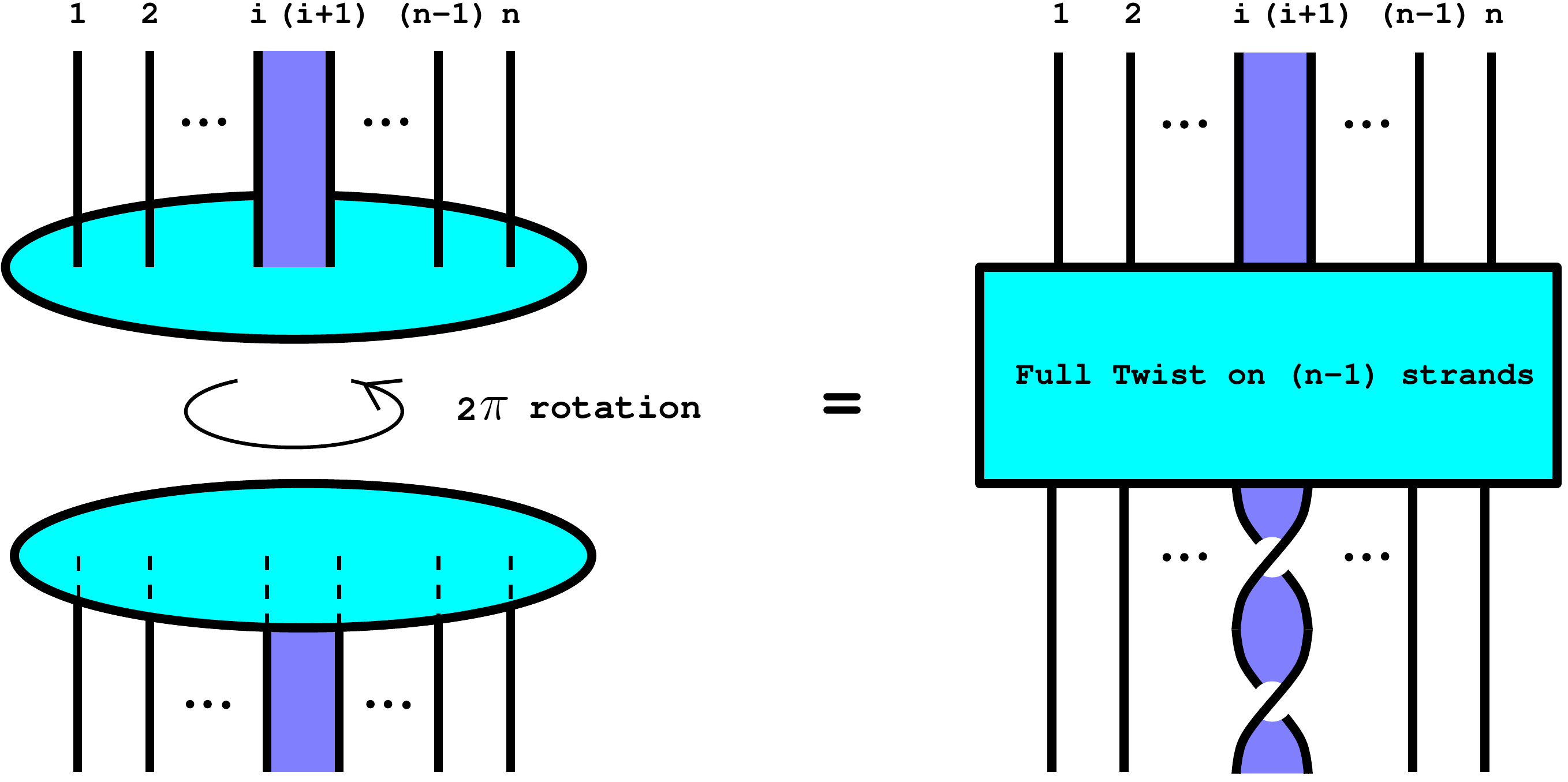}
\caption{Full twist on $n$ strands, with two consecutive strands grouped together as one strand.}
\label{basecaseone}
 \end{center}
 \end{figure}

\begin{figure}
\begin{center}
\includegraphics[height=2in]{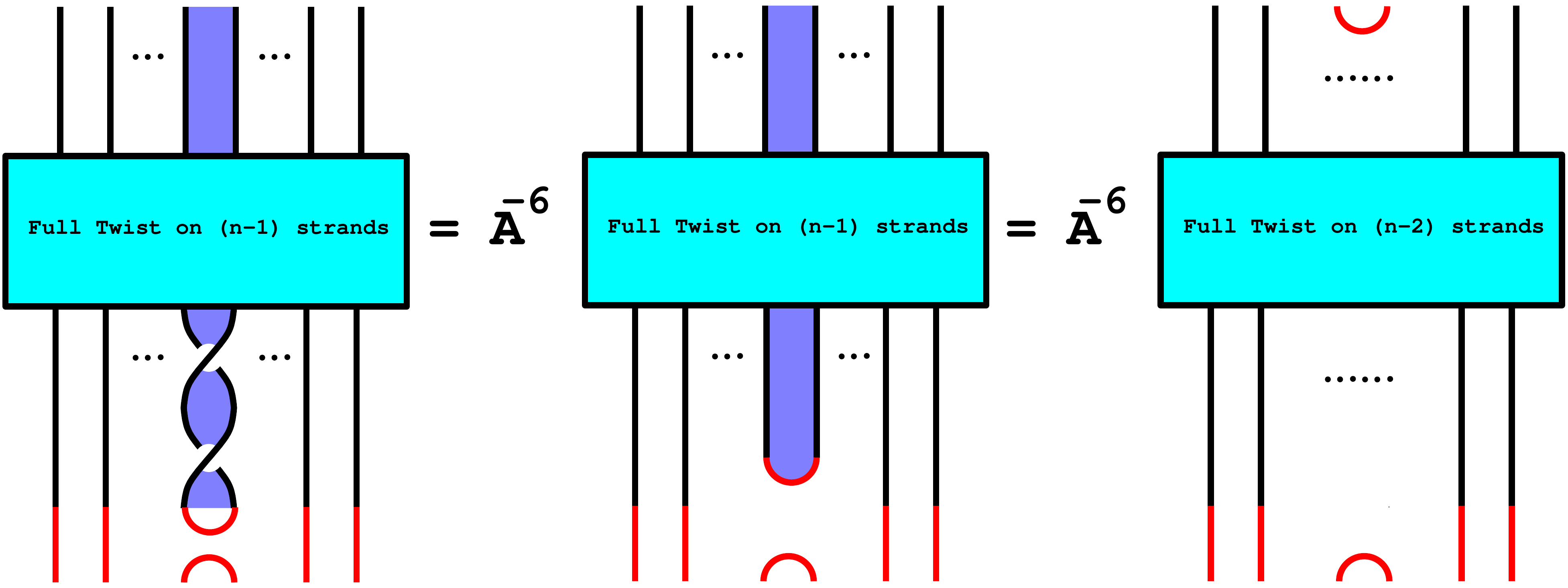}
\caption{$\D_n e_i = A^{-6} E_i$ }
\label{basecasetwo}
 \end{center}
 \end{figure}

 We will prove the claim by induction on the length of $h$ as a
 product of $e_i$'s.  For the base case, $h = e_i$. If $E_i$ denotes
 $e_i$ with its $n-2$ through strands given a full right twist, then
 from Figures \ref{basecaseone} and \ref{basecasetwo} we see that $\D_n
 e_i = A^{-6} E_i$.

Assuming the claim holds for $h'=e_{i_1} \ldots e_{i_{r-1}}$, we must
show that it holds for $h=h' e_{i_r}$.  For the standard $TL_n$ basis,
we have that $i_r$ is distinct from $i_1,\ldots,i_{r-1}$, so that
${\rm cap}(e_{i_r}) \nsubseteq {\rm cap}(h') $. This gives us two
cases (see Figure \ref{fulltwisteval}):

\begin{figure}
\begin{center}
\includegraphics[height=1.2in]{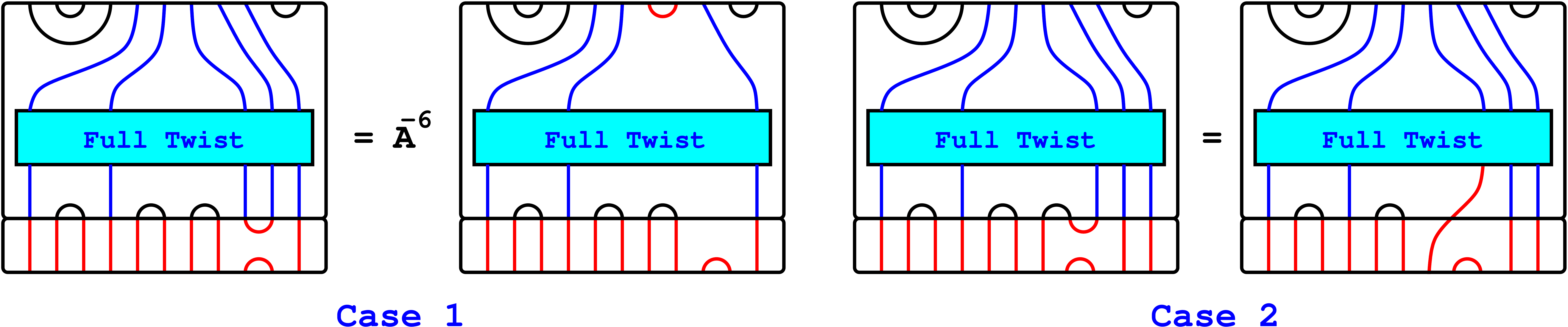}
\caption{Two cases for Lemma \ref{fulltwistevaluation}}
\label{fulltwisteval}
 \end{center}
 \end{figure}

 {\bf Case 1:} ${\rm cap}(e_{i_r}) \subseteq {\rm feet}(h')$. In this
 case, $|{\rm cap}(h)|=|{\rm cap}(h')|+2$.  The full twist on $|{\rm
   feet}(h')|$ strands is paired with $e_{i_r}$ and results in $A^{-6}
 \D_{|{\rm feet}(h')| -2}$. Hence, $\D_n h= A^{-6(k-1)}H'e_{i_r}=A^{-6k} H$.

 {\bf Case 2:} ${\rm cap}(e_{i_r}) \nsubseteq {\rm feet}(h')$. In this
 case, ${\rm cap}(e_{i_r})$ is split between ${\rm feet}(h')$ and
 ${\rm cap}(h')$.  The number of caps and the number of through
 strands of $h$ are both the same as those in $h'$.
Hence, $\D_n h= A^{-6k}H'e_{i_r}=A^{-6k} H$.
\done

\begin{lemma}
\label{maxpower}
The highest power of $A$ in $\kb{\cl{\D_n h_i}}$ for any $i>0$ is less than or equal to $n^2-3n-4$. 
\end{lemma}
\pf Following the notation in Lemma \ref{fulltwistevaluation}, $k =
\#$cups in $h_i$, so that $1\leq k \leq [n/2]$.  By Lemma
\ref{fulltwistevaluation}, $\D_n h_i = A^{-6k} H_i$ if $i>0$, where
$H_i$ also has $k$ cups.  The closure of $H_i$ will result in $k'$
cups paired with caps to produce loops, and $k''$ cups pulled
through the full twist with a factor of $A^{-6k''}$, where $0\leq
k,k''\leq k$.  So $\cl{H_i}$ will have a full twist on $m'$
strands, where $0\leq m'\leq m$. 
Thus, $\kb{\cl{\D_n h_i}} = A^{-6k}\kb{\cl{H_i}}=A^{-6(k+k'')}\delta^{k'}\kb{\cl{\D_{m'}}}$.

The highest power of $A$ is $(-6k -6k'' + 2k' +
(m')^2 + m'-2)$, which is maximized when $k''=0,\ k'=k$ and
$m'=m$.  Since $2k+m=n$, we have $m^2+m-2-4k = (n-2k)^2+(n-2k)-2-4k$.
The function $f(k)=(n-2k)^2+(n-2k)-2-4k$ has an absolute minimum at
$k= \frac{n}{2}+\frac{3}{4}$ and is decreasing for $k<
\frac{n}{2}+\frac{3}{4}$. Since $1 \leq k \leq [n/2]<
\frac{n}{2}+\frac{3}{4}$, $f(k) \leq f(1) = (n-2)^2 +(n-2) -2 -4=
n^2-3n -4$.  \done

\section{Proof of Theorem \ref{mainthm2}}

% We now extend our results to the colored Jones polynomials. 
For any link diagram $D$, let $D^{(r)}$ denote its {\em blackboard
  framed} $r$--cable.  Let $D_n$ denote the standard diagram of the
closure of a full positive twist on $n$ strands with a positive kink
on each strand.

\begin{lemma}
\label{cabling_fulltwist}
 $\displaystyle D_n^{(r)}=D_{nr}$.
\end{lemma}
\pf The {\em belt trick} (see, e.g., \S 2.2-2.4 of \cite{Wenzl_3mfld}) implies:
\begin{center}
\includegraphics[height=1.25in]{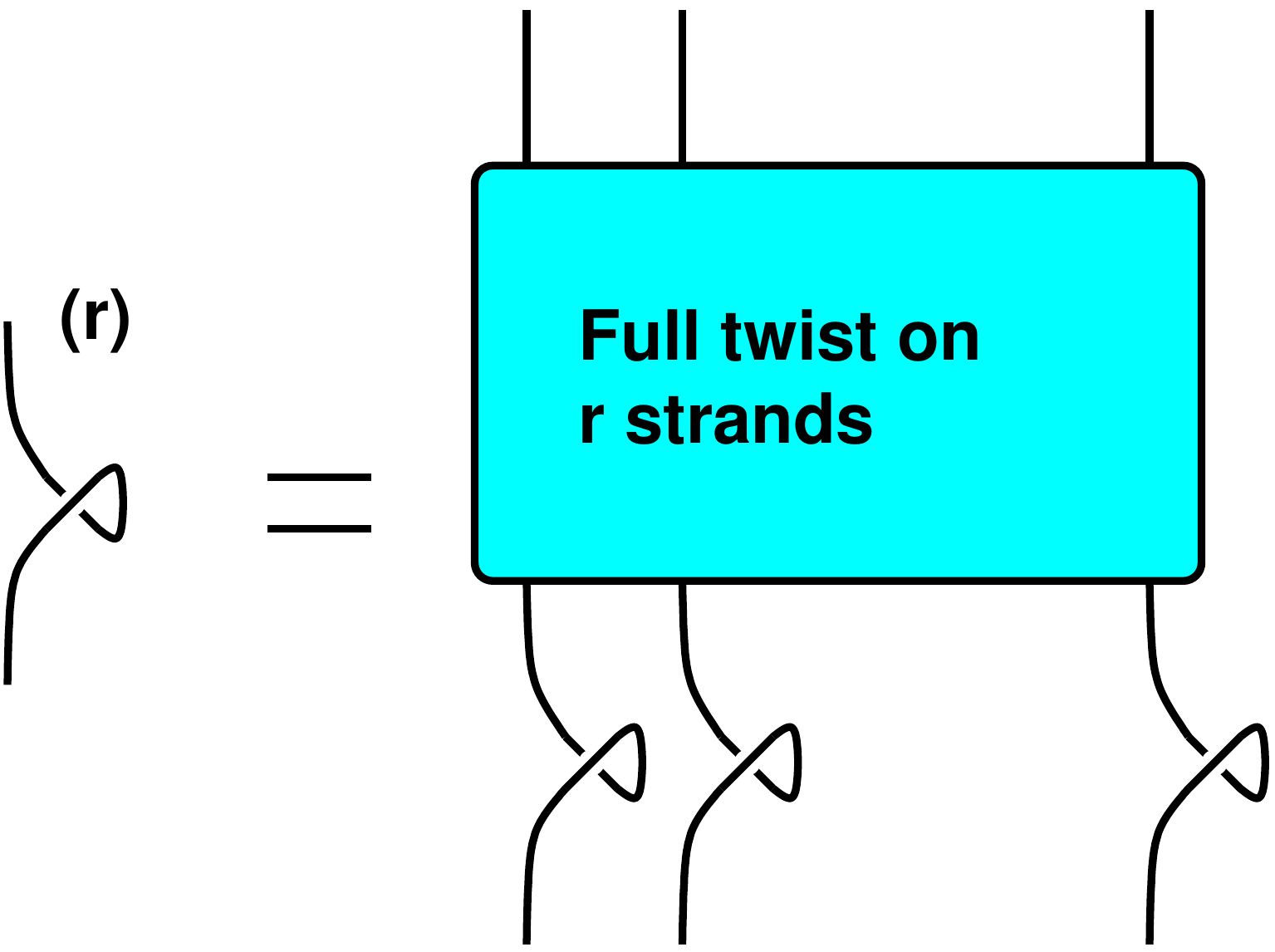}
\end{center}
Thus, $D_n$ is planar isotopic to
$(\includegraphics[height=.115in]{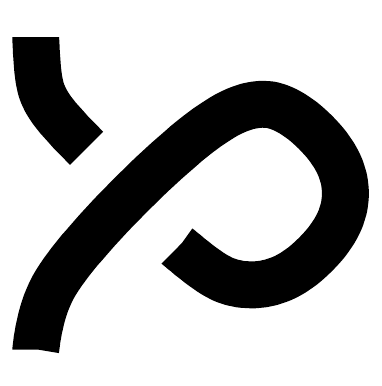})^{(n)}$. Hence
$\displaystyle D_n^{(r)}= ((\includegraphics[height=.115in]{pkink})^{(n)})^{(r)} =
(\includegraphics[height=.115in]{pkink})^{(nr)}=D_{nr}.$
\done

For $N\geq 0$, let $S_N(x)$ be the Chebyshev polynomials defined by
$$ \displaystyle{S_N(x)=\sum^{[N/2]}_{j=0}(-1)^j\binom{N-j}{j}x^{N-2j}}. $$
For a link diagram $D$, let $S_N(D)$ be a linear combination of
blackboard cablings of $D$.
We can define the colored Jones polynomial, as in equation (\ref{colJP}), by
a related expression in terms of the blackboard framing for $D$ 
(see, e.g., \cite{DL2}):
\begin{equation}
\label{coloredjp}
J_{N+1}(D;A^{-4})=\left( (-1)^N A^{N^2+2N}\right)^{-w(D)}(-1)^{N-1} (A^2+A^{-2}) \kb{S_N(D)}
\end{equation}

As before, let $\beta'$ be a positive $n$--braid with $c$ crossings,
$\beta=\D_n \,\beta'$, and $L=\cl{\beta}$.  Let $D$ be a diagram of
$L$ given by the closure of $\beta$ with a positive kink on every
strand after the full twist. Note that the writhe of $D$,
$w(D)=n^2+c$. By Lemma \ref{cabling_fulltwist}, $D^{(r)}$ is the
closure of the braid $\D_{nr}\beta'^{(r)}$ with a positive kink on
every strand following the full twist.

Let $r(A)=(-1)^{N(w(D)+1)+1} A^{-(N^2+2N)w(D)} (A^2+A^{-2})$.  By equation (\ref{coloredjp}),
\begin{align}
\label{coloredjp2}
& (1-A^{-8})J_{N+1}(L;A^{-4}) = (1-A^{-8})\, r(A)
 \sum^{[N/2]}_{j=0}(-1)^j\binom{N-j}{j} \kb{D^{(N-2j)}}\qquad \notag\\
&\qquad = r(A)  \sum^{[N/2]}_{j=0}(-1)^j\binom{N-j}{j} (-A^3)^{n(N-2j)}
(1-A^{-8}) \kb{\overline{\D_{n(N-2j)}\beta'^{(N-2j)}}}
\end{align}

By equation (\ref{bracketofbeta}) with $\ell=0$,
\begin{align}
\label{bracketofbeta1}
(1-A^{-8})\kb{\overline{\D_n \beta'}}&= 
(-1)^{n+1}\left( q_3(A)+(n-2)A^{n^2-3n-6+c} + A^{n^2+n-2+c} \right)\notag \\
&=(-1)^{n+1} \left( q_4(A)A^{n^2-3n-6+c} + A^{n^2+n-2+c}\right)
\end{align}
where $q_3(A)$ and $q_4(A)$ are Laurent polynomials such that
$\max\deg(q_3(A))\leq n^2-3n-6+c$ and $\max\deg(q_4(A))\leq 0$.

Let $\displaystyle \ d_1(j)=(N-2j)^2(n^2+c)-6 \quad {\rm and} \quad d_2(j)=(N-2j)^2(n^2+c)+4n(N-2j)-2$.

Note that $d_1(j)$ and $d_2(j)$ are both quadratic functions of $j$.
As $j$ increases from $0$ to $[N/2]$, they decrease and $d_2(j) > d_1(j)$.

\begin{lemma}
\label{dlemma}
$d_1(j)> d_2(j+1)$ for $0 \leq j \leq [N/2]-1$. 
\end{lemma}
\pf Note that $n \geq 3$ and $N \geq 3$. 
\begin{align*} 
  d_1(j)-d_2(j+1) &=\left((N-2j)^2(n^2+c)-6\right) -\left((N-2j-2)^2(n^2+c)+4n(N-2j-2) -2\right) \\
  &=4(N-2j)(n^2+c)-4(n^2+c)-4n(N-2j-2)-4 \\
  &=4n(n(N-2j)-(N-2j)-n)+4c(N-2j-1)+8n-4\\
  &=4n(n-1)(N-2j-1)+4c(N-2j-1)+4(n-1)\\
  & > 0 {\rm\ for\ } 0\leq j \leq [N/2]-1
\end{align*}
\done

Let $b_j=(-1)^j\binom{N-j}{j}$.  Continuing from equation (\ref{coloredjp2}), we now have
\begin{align*}
(1-A^{-8})J_{N+1}(L;A^{-4})
& = r(A)  \sum^{[N/2]}_{j=0}b_j (-A^3)^{n(N-2j)}
(1-A^{-8}) \kb{\overline{\D_{n(N-2j)}\beta'^{(N-2j)}}} \\
&= (-1)^{nN}r(A)  \sum^{[N/2]}_{j=0} (-1)^{n(N-2j)+1} b_j
\big[ \bar p_{N-2j}(A)A^{d_1(j)} + A^{d_2(j)} \big] \ \ {\rm \ using\ (\ref{bracketofbeta1})} \\
& = -r(A)  \sum^{[N/2]}_{j=0} b_j
\big[  A^{d_1(j)} \bar p_{N-2j}(A) + A^{d_2(j)}   \big] \\
&=-r(A) \big[ A^{d_1(0)} p_{N}(A) + A^{d_2(0)} \big] {\rm \ using\ Lemma\ \ref{dlemma} \ } \\
\end{align*}
where $p_N(A)$ and $\bar p_{N-2j}(A)$ are Laurent polynomials 
with $\max\deg p_N(A) \leq 0$ and \\
 $\max\deg(\bar p_{N-2j}(A)) \leq 0$ for $0 \leq j \leq [N/2]$, respectively. 

Substituting $t=A^{-4}$ and, in $r(A)$, $w(D)=n^2+c$,
\begin{align*}
(1-t^2)J_{N+1}(L;\,t) &= (-1)^{-N(n^2+c-1)}t^{\frac{(N^2+2N)(n^2+c)}{4}} (t^{1/2}+t^{-1/2})\left(t^{-d_2(0)/4} + t^{-d_1(0)/4}p_N(t)\right) \\
 &= (-1)^{-N(n^2+c-1)}(t^{1/2}+t^{-1/2})
\left(t^{\frac{N(n^2+c-2n)+1}{2}} + t^{\frac{N(n^2+c)+3}{2}}p_N(t)\right) \\
(1-t)J_{N+1}(L;\,t)&=(-1)^{N(n^2+c-1)}t^{\frac{N(n^2+c-2n)}{2}} 
\left(1+t^{nN+1}p_N(t)\right)
\end{align*}
where $p_N(t)$ is a polynomial. 
This completes the proof of the first statement of Theorem \ref{mainthm2}.  
Now, the second statement follows from
$$J'_N(L;\,t)=\frac{J_N(L;\,t)}{[N]}=\frac{t^{\frac{N-1}{2}}(1-t)J_N(L;\,t)}{(1-t^{N})}.$$

To obtain the tail (without denominators), we suppress $\cl{\beta}$ in the notation, 
\begin{align*}
\frac{1+t^{nN+1}p_N(t)}{1-t} %& = \frac{1-t^{nN+1} + t^{nN+1}(1+p_N(t))}{1-t} \\
&= \frac{1-t^{nN+1}}{1-t} + \frac{t^{nN+1}(1+p_N(t))}{1-t} \\
&= \sum_{i=0}^{nN} t^i + t^{nN+1}q_N(t) \\ \\
\frac{1+t^{nN+1}p_N(t) }{1-t^{N+1}} &= \frac{1-t^{n(N+1)}}{1-t^{N+1}} 
+ \frac{t^{nN+1}(t^{n-1}+p_N(t) )}{1-t^{N+1}}\\
&= \sum_{i=0}^{n-1}t^{i(N+1)}  + t^{nN+1} q_N(t) \\
\end{align*}
This completes the proof of Theorem \ref{mainthm2}.
\done  

\bibliography{references}
\bibliographystyle{plain}

\end{document}